# Tight Chromatic Upper Bound for $\{3K_1, K_1+C_4\}$-free Graphs

Medha Dhurandhar
mdhurandhar@gmail.com

**Abstract**

Problem of finding an optimal upper bound for a chromatic no. of $3K_1$-free graphs is still open and pretty hard. It was proved by Choudum et al that upper bound on the chromatic no. of $\{3K_1, K_1+C_4\}$-free graphs is $2\omega$. We improve this by proving that if G is $\{3K_1, K_1+C_4\}$-free, then $\chi \leq \frac{3\omega}{2}$ where $\omega$ is the size of a maximum clique in G. Also we give examples to show that the bound is tight.

**1. Introduction:**

In [1], [2], [4], [5] chromatic bounds for graphs are considered especially in relation with $\omega$ and $\Delta$. Gyárfás [6] and Kim [7] show that the optimal $\chi$-binding function for the class of $3K_1$-free graphs has order $\omega^2/\log(\omega)$. If we forbid additional induced subgraphs, the order of the optimal $\chi$-binding function drops below $\omega^2/\log(\omega)$. We consider in this paper only finite, simple, connected, undirected graphs. The vertex set of G is denoted by V(G), the edge set by E(G), the maximum degree of vertices in G by $\Delta(G)$, the maximum clique size by $\omega(G)$ and the chromatic number by $\chi(G)$. N(u) denotes the neighbourhood of u and $\overline{N(u)} = N(u) + u$. In [1] it was proved that is $2\omega$ is an upper bound on the chromatic no. of $\{3K_1, K_1+C_4\}$-free graphs and the problem of finding tight chromatic upper bound for a $\{3K_1, K_1+C_4\}$-free graph was stated as open.

For further notation please refer to Harary[3].

**2. Main Result:**

The following Lemma defines the structure of $\{3K_1, K_1+C_4\}$-free graphs.

**Lemma 1:** If G is $\{3K_1, K_1+C_4\}$-free, then $G = \bigcup_1^j M_i$ where $1 \leq j \leq 4$ and $\langle M_i \rangle$ is complete $\forall\ i$, $1 \leq i \leq j$ with $|M_1| \leq \omega$, $|M_2| \leq \omega$, $|M_3| \leq \omega-1$, and $|M_4| \leq \omega-1$.

Proof: If G is complete, then $G = M_1$. Next let G be not complete and v, w $\in$ V(G) be s.t. vw $\notin$ E(G). Let $A = \{x \in V(G)/xv, xw \in E(G)\}$, $B = \{x \in V(G)/xv \in E(G), xw \notin E(G)\}$, $C = \{x \in V(G)/xw \in E(G), xv \notin E(G)\}$. Let $A_1 \subseteq A$ be s.t. $\langle A_1 \rangle$ is a maximal clique in $\langle A \rangle$. Let $A_2 = A-A_1$. As G is $3K_1$-free, $\langle B \rangle$ and $\langle C \rangle$ are complete. Now $\langle A_2 \rangle$ is complete. Else let $a_{2i} \in A_2$ (i=1, 2) be s.t. $a_{21}a_{22} \notin E(G)$. Clearly as G is $3K_1$-free, $|A_1| > 1$ and $\exists\ a_{1i} \in A_1$ (i=1, 2) s.t. $a_{1i}a_{2i} \notin E(G)$ for i = 1, 2. But then $\langle a_{12}, a_{11}, v, a_{21}, w \rangle = K_1+C_4$, a contradiction. Hence $\langle A_2 \rangle$ is complete. Let $M_1 = A_1 \cup v$, $M_2 = A_2 \cup w$, $M_3 = B$ and $M_4 = C$.

This proves the Lemma.

**Main Result:** If G is $\{3K_1, K_1+C_4\}$-free, then $\chi \leq \frac{3\omega}{2}$.

Proof: Let if possible G be a smallest $\{3K_1, K_1+C_4\}$-free graph with $\chi > \frac{3\omega}{2}$. Let $v \in V(G)$ be s.t. deg $v = \Delta$. If $\Delta = |V(G)|-1$, then $\omega(G) = \omega(G-v)+1$ and by minimality $\chi(G-v) \leq \frac{3\omega(G-v)}{2}$. Thus $\chi(G) \leq \chi(G-v)+1 \leq \frac{3\omega(G-v)}{2}+1 < \frac{3\omega}{2}$, a contradiction. Hence deg $v = \Delta < |V(G)|-1$. Let $w \in V(G)$ be s.t.

$vw \notin E(G)$. Let $A = \{a \in V(G) / av, aw \in E(G)\}$, $B = \{b \in V(G) / av \in E(G)$ and $aw \notin E(G)\}$, and $C = \{c \in V(G) / cv \notin E(G)$ and $aw \in E(G)\}$. As $G$ is $3K_1$-free, $V(G) = A \cup B \cup C$. Let $A_1 \subseteq A$ be s.t. $<A_1>$ is a maximal clique in $<A>$. Let $A_2 = A - A_1$. Let $A_j = \bigcup_{i=1}^{|A_i|} A_{ji}$ for $j = 1, 2$. Now $<A_2>$ is complete as shown in **Lemma 1**. Again $a_{1i}a_{2j} \notin E(G)$ $\forall$ j (else let $a_{1i}a_{2j} \in E(G)$ and $a_{1k}a_{2j} \notin E(G)$. Then $<a_{1i}, a_{2j}, v, a_{1k}, w> = K_1 + C_4$). Next let $B_{1k} = \{b \in B / ba_{mn} \notin E(G)$ for some n$\}$ where k, m $\in \{1, 2\}$, k $\neq$ m, and $B_2 = B - \bigcup_{i=1}^{2} B_{1i}$. Then $<B_{11} \cup A_1 \cup B_2>$, $<B_{12} \cup A_2 \cup B_2>$ are complete and $|A_1| + |B_{11}| + |B_2| + 1 \leq \omega$, $|A_2| + |B_{12}| + |B_2| + 1 \leq \omega$. Thus $\Delta = |A_1| + |B_{11}| + |B_2| + |A_2| + |B_{12}| \leq \omega - 1 + \omega - 1 - |B_2|$. Hence $\chi \leq \lceil \frac{\Delta + \omega + 1}{2} \rceil \leq \frac{3\omega}{2}$, a contradiction.

This proves the result.

**Remarks:** Examples of $\{3K_1, K_1+C_4\}$-free graphs where the chromatic no. equals the upper bound.
1. Let $\mathbf{G} = \sum_{i=1}^{2k} G_i$ where $G_i = C_5 \forall$ i. Then $\omega = 2k$ and $\chi = 3k = \frac{3\omega}{2}$.
2. Let $\mathbf{G} = \sum_{i=1}^{2k} G_i$ where $G_i = W_6 \forall$ i. Then $\omega = 3k+1$,

**References:**
[1] "Linear Chromatic Bounds for a Subfamily of 3K1-free Graphs", S. A. Choudum, T. Karthick, M. A. Shalu, Graphs and Combinatorics 24:413–428, 2008
[2] "On the divisibility of graphs", Chinh T. Hoang, Colin McDiarmid, Discrete Mathematics 242, 145–156, 2002
[3] Harary, Frank, *Graph Theory* (1969), Addison–Wesley, Reading, MA.
[4] "ω, Δ, and χ", B.A. Reed, J. Graph Theory 27, pp. 177-212, 1998
[5] "Some results on Reed's Conjecture about ω, Δ and χ with respect to α", Anja Kohl, Ingo Schiermeyer, Discrete Mathematics 310, pp. 1429-1438, 2010